\documentclass[11pt,reqno]{article}
\usepackage{graphicx,epsfig}
\usepackage{amsmath}

\title {$\omega$-Periodic graphs}

\author{ Itai Benjamini \and Chris Hoffman}

\usepackage{latexsym}
\usepackage{amsfonts}
\usepackage{amsmath}
\usepackage{amsthm}
\usepackage{graphicx}

\newtheorem{thm}{Theorem}
\newtheorem{defn}[thm]{Definition}
\newtheorem{lemma}{Lemma}

\newtheorem{prop}[thm]{Proposition}

\newcommand{\be}{\begin{equation}}
\newcommand{\ee}{\end{equation}}

\newcommand{\Z}{{\mathbb Z}}

\begin{document}

\maketitle

\begin{abstract}
$\omega$-periodic graphs are introduced and studied.  These are
graphs which arise as the limits of periodic extensions of the
nearest neighbor graph on the integers. We observe that all
bounded degree $\omega$-periodic graphs are ameanable. We also
provide examples of $\omega$-periodic graphs which have
exponential volume growth, non-linear polynomial volume growth and
intermediate volume growth.
\end{abstract}

\section{Introduction}

In \cite{M} Milnor asked the following question.  Does there exist
a finitely generated group $G$ such that the volume of the ball of
radius $n$ about the identity of Cayley graph of $G$ grows faster
than polynomially but slower than exponentially? This question was
answered by Grigorchuk, who constructed a family of groups whose
Cayley graphs have intermediate growth \cite{G}. See \cite{DL} for
a nice description of these groups.

Graphs that have intermediate volume growth also have a strong
connection with long range percolation models in probability. In
long range percolation on $\Z$, a random graph is constructed with
$\Z$ as its vertex set. The measure is determined by  a sequence
$p_n$. For each pair $(u,v) \in \Z \times \Z$ there exists an edge
$e_{u,v}$ between $u$ and $v$ with probability $p_{|u-v|}$. The
existence of an edge between $u$ and $v$ is determined only by the
distance between $u$ and $v$ and is independent of the existence
of edges between any other pairs of vertices. Long range
percolation on $\Z$ was introduced and studied in \cite{s},
\cite{NS} and \cite{AN} and is commonly used as a model for social
networks.

Given a sequence $p_n$ these papers studied the probability that
an infinite connected subgraph exists. The more recent papers
\cite{Be}, \cite{BB}, \cite{Ga} and \cite{Bi} considered the case
when there is a unique infinite connected subgraph a.s.\ and
studied the volume growth of this graph. In particular they
considered the case that $p_1=1$ and $p_{n} = \beta n^{\alpha}$
for $n>1$. For these sequences when $\alpha >2$ the random graph
has linear volume growth a.s. When $2> \alpha >1$, the random
graph has intermediate growth a.s., yet large intervals admit
polynomially small boundaries. And when $ \alpha \leq 1$, all
degrees are infinite. It is conjectured in \cite{BB} that when
$\alpha =2$ one gets polynomial volume growth with powers
depending on $\beta$.

In this paper we present a simple and natural graph, $G$, that has
intermediate volume growth.  Our graph is in some sense a hybrid
of the Cayley graphs of the Grigorchuk groups and the graphs from
long range percolation.  Like the graphs generated by long range
percolation our graph is constructed as an extension of the
nearest neighbor graph on $\Z$.  However its description is
deterministic and, like Cayley graphs, has much more regularity
than those random graphs. It does not have the full symmetry that
Cayley graphs posses, it has small bottlenecks and in particular
it is not transitive.

In addition to the study of one particular graph we also consider
a broad family of graphs that contains $G$.  This is the set of
all graphs which are constructed as limits of periodic graphs on
$\Z$. We call these  $\omega$-periodic graphs.  We study some of
the  coarse geometric properties shared by all $\omega$-periodic
graphs and then give a few more examples that illustrate the
possible volume growth of $\omega$-periodic graphs. In particular
we show that there are $\omega$-periodic graphs of (non-linear)
polynomial growth as well as ones with exponential volume growth.

\section{The Basic Example}
\label{inter}

The vertices of $G$ are the integers, $\Z$. Define the sets of
edges
$$E_0 =\{(i,i+1) : i \in \Z\}$$
and
$$
E_k=\{\left(2^k (n -1/2), 2^k(n+1/2)\right)\},
$$
for all $n \in \Z$ and $k >0$.

The graph $G$ has edges $E= \cup_{k\geq 0} E_k.$ $G$ is
$\omega$-periodic because it is the union of $G_K$ which has edges
$\cup_0^k  E_k$.  We refer to the edges in $E_k$ as the $k$th
layer. The degree of every vertex of $G$ (except for 0) is $4$.
More useful than our description is the picture below.

\newpage

\begin{figure}
\centerline{ \scalebox{.6}{\includegraphics[width=20cm]{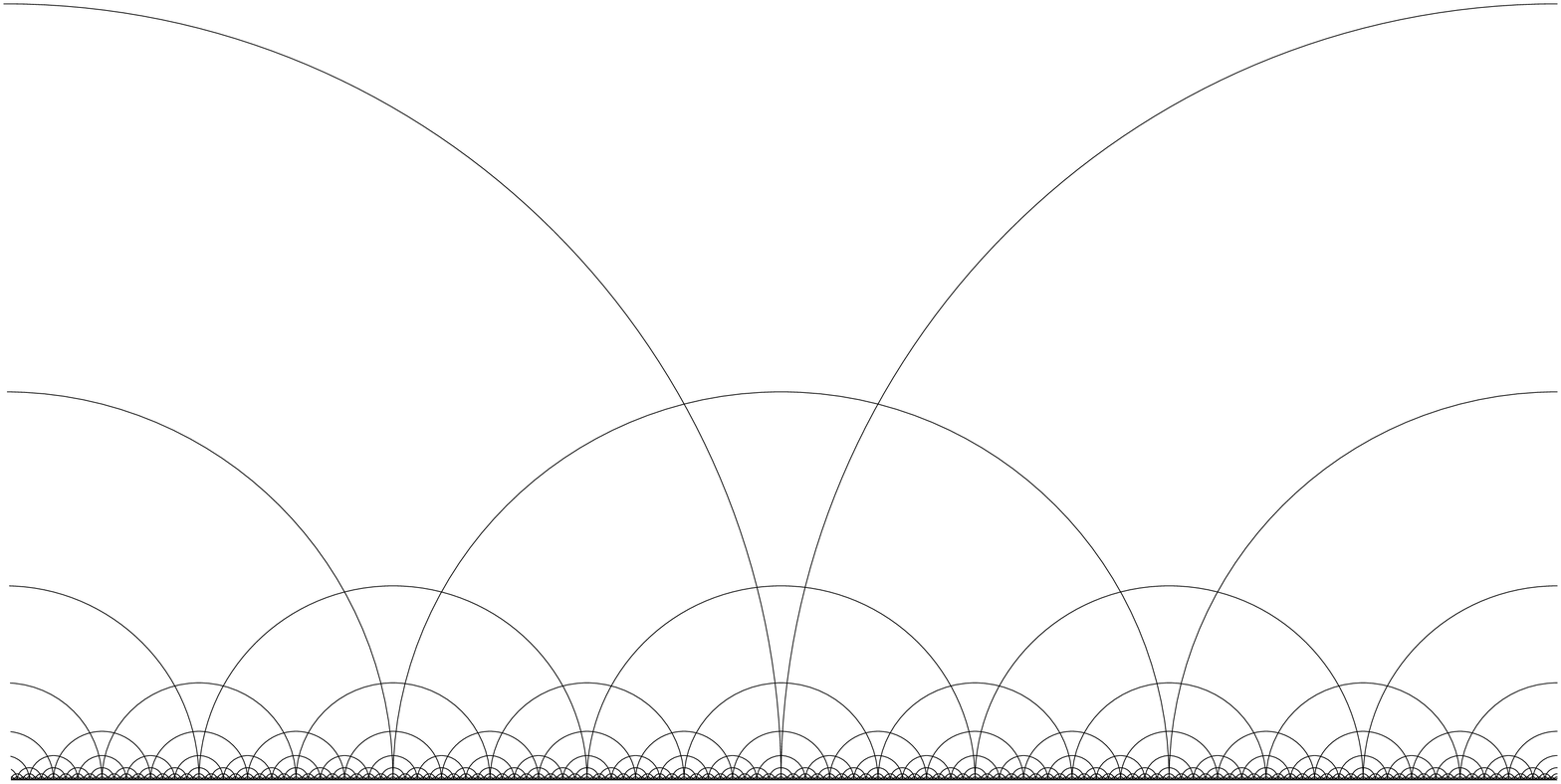}}}
\caption{The basic example}
\end{figure}

In order to calculate the volume growth of $G$ we would like to
calculate the length minimal path between any two integers and
then use that to estimate the volume growth.  To find a minimal
path between two integers $u$ and $v$ we take the following
approach. We pick an integer $k$ we move as quickly as possible
from $u$ to a vertex in the $k$th layer. Then we move in the $k$th
layer and finally we go from the $k$th layer to $v$. It is easy to
show that a minimal path must take such an approach.

The problem is given $u$ and $v$ how do we identify the optimal
$k$.  In general we do not know how to answer that question but we
are able to calculate the length minimal paths between 0 and
points of the form $2^n$.  (This is called $x_{n+1}$.)  Our main
tools are induction and the symmetries of the graphs $G_k$. By
knowing the distance from 0 to $2^k$ for all $k<n$ we can
determine the distance from 0 to $2^n$.  The inductive
relationship is given in Lemma \ref{one} while the formula is
determined explicitly in Lemma \ref{two}.

Using this information along with the symmetries of $G_k$ we are
able to determine the growth rate of $|B_j(0)|$, the number of
vertices within distance $j$ of 0. Although $|B_j(0)|$ does not
have a simple formula we show in Lemma \ref{growth} that $|B_j(0)|
\approx j^{.5 \log j}$ and determine $|B_j(0)|$ to within a factor
of $16j^2$ .

To analyze the growth rate of the ball centered at 0 we make the
following definitions. For $i\geq 1$ and $j \in \Z$ let $x_{i,j}$
be the distance from 0 to $(j+1/2)2^i$ in the graph $G_i$. The
next lemma gives us an inductive relationship for $x_i$ based on
$x_j$, $j<i$.  In Lemma \ref{two} we will calculated $x_i$
explicitly.
\begin{lemma}
\label{one}
$$x_{i,j}=x_{i,0}+|j+1/2|-1/2.$$
For all $i>2$
$$x_i=x_{i,0}=\min_{0< k < i}2x_k+2^{i-k-1}-1.$$

\end{lemma}

\begin{proof}
The proof is by induction on $i$. It is easy to check for $i=1,2$
that the first formula is true.  Fix $i$ and $j$. Assume the lemma
is true for all $k<i$ and for all $j$.

Fix $j$.
Let $P$ be an oriented path in
$G_i$ from 0 to $(j+1/2)2^i$ which has minimal length. Let $k$
be the largest integer
such that an edge of $P$ is in $E_k$. It causes no loss of generality to
assume that $k<i$. (This is because there is a first point of the form
$(j'+1/2)2^i$ in $P$.  If the lemma is not true for $i$ and $j$ then it
is also not true for $i$ and $j'$.)

Since $P$ has minimal length
and $i>2$ then it is clear that $k>0$. Divide $P$ up into three
parts, $P_1,P_2,$ and $P_3$, as follows. Let $n_1$ be the first
point in $P$ of the form $(n_1+1/2)2^k$. Let $n_2$ be the last
point in $P$ of the form $(n_2+1/2)2^k$.  Then $P_1$ is the
portion of $P$ connecting 0 to $(n_1+1/2)2^k$, $P_2$ connects
$(n_1+1/2)2^k$ to $(n_2+1/2)2^k$, and $P_3$ connects
$(n_2+1/2)2^k$ to $(j+1/2)2^i$. Then we have that
\begin{eqnarray*}
|P| &=& |P_1|+|P_2|+|P_3|\\
    &\geq& x_{k,n_1}+(n_2-n_1)+x_{k,2^{i-k-1}-n_2-1+j2^{i-k}}\\
    &\geq& x_k + |n_1-1/2|+1/2+(n_2-n_1)+x_k+
           |2^{i-k-1}-n_2+j2^{i-k}-1/2|-1/2\\
    &\geq& 2x_k + n_1+(n_2-n_1)+ 2^{i-k-1}-n_2-1+j2^{i-k}\\
    &\geq& 2x_k + 2^{i-k-1}-1+j2^{i-k}\\
    &\geq& \min_{0<k<i} 2x_k + 2^{i-k-1}-1+j2^{i-k}.
\end{eqnarray*}
The existence of a path of the minimum distance is easy to
construct.
\end{proof}

We now calculate $x_i$ exactly.
Let
$$y_{n}=\frac{n^2+3n+2}{2}$$
and
$$z_n=n2^{n}+1.$$

\begin{lemma}
\label{two}
For all $i$ and $n$, if $y_n<i \leq y_{n+1}$ then
$$x_i=z_{n+1}-(y_{n+1}-i)2^n.$$  In particular
$x_{y_n}=z_n$.
\end{lemma}

\begin{proof}
The proof is by induction.
It is easy to check that the lemma is true for all $i\leq y_1=3$.

Now assume that the lemma is true for all $j<i$.
Note that this implies that the sequence
$$x_j-x_{j-1}$$
is nondecreasing for all $j$, $2 \leq j \leq i-1$.

Let
$$f(k)=f(i,k)=2x_k + 2^{i-k-1}-1.$$
By Lemma \ref{one} $x_i= \min_{k<i} f(k).$
Let $n$ be the largest integer such that $y_n<i$.  We break the proof up
into two cases.  The first is when $i<y_{n+1}$ and the second is when
$i=y_{n+1}$.

\noindent
{\bf Case 1:}
We show that the minimum of $f(k)$ occurs at two values, $i-n-1$ and
$i-n-2$.  More specifically we show that $f(k)$ is decreasing up to $i-n-1$
 and increasing afterwards.  Let $m=i-n-1$.
Since $y_n<i<y_{n+1}$ we have that $y_{n-1}<m\leq y_n$. Thus
$$x_{m}-x_{m-1}=2^{n-1}$$
and
$$x_{m+1}-x_{m} \geq 2^{n-1}.$$

We now calculate for $j<m$
\begin{eqnarray*}
f(j)-f(j-1)
&=    &2(x_{j}-x_{j-1})+(2^{i-j-1}-2^{i-j})\\
&\leq &2(x_{m}-x_{m-1})-2^{i-j-1}\\
&\leq &2 \cdot 2^{n-1}-2^{i-m}\\
&\leq &2^n-2^{n+1}\\
&<    &0.
\end{eqnarray*}

We also have that
\begin{eqnarray*}
f(m)-f(m-1)
&=&2(x_m-x_{m-1})+(2^n-2^{n+1})\\
&=&2 \cdot 2^{n-1}-2^n\\
&=&0
\end{eqnarray*}
For $l>m$ we have

\begin{eqnarray*}
f(l+1)-f(l)
&=    &2(x_{l+1}-x_{l})+(2^{i-l-2}-2^{i-l-1})\\
&\geq &2(x_{m+1}-x_{m})-2^{i-l-2}\\
&\geq &2(2^{n-1})-2^{i-m-2}\\
&\geq &2^{n}-2^{n-1}\\
&>&0.
\end{eqnarray*}
From these three calculations
it is clear that $f$ obtains its minimum at $m-1$
and $m$. It is easy to check that the induction hypothesis gives
the right value for $x_i$.

\noindent
{\bf Case 2:}
Now we have that $i=y_{n+1}$. We claim that in this case the unique
minimum of $f$ occurs at
$m=y_{n}$.
By the induction hypothesis we
have that
$$x_{m}-x_{m-1}=2^{n-1}$$
and
$$x_{m+1}-x_{m} = 2^{n}.$$

We now calculate for $j<m$
\begin{eqnarray*}
f(j)-f(j-1)
&=&      2(x_j-x_{j-1})+(2^{i-j-1}-2^{i-j})\\
&\leq &  2(x_m-x_{m+1})-2^{i-j-1}\\
&=    &  2 \cdot 2^{n-1}-2^{i-m-1}\\
&=    &  2^{n}-2^{n+1}\\
&=    &  2^n>0.
\end{eqnarray*}

For $l>m$
\begin{eqnarray*}
f(l+1)-f(l)
&=&       2(x_{l+1}-x_l)+2^{i-l-2}-2^{i-l-1}\\
& \geq &  2(x_{m+1}-x_{m})-2^{i-l-2})\\
&\geq  &  2^{n+1}-2^{i-m-2}\\
&\geq  &  2^{n+1}-2^{n}\\
&=&       2^n\\
&>&0.
\end{eqnarray*}
From these calculations it is clear that $f$ obtains its minimum at $m$.

Thus
\begin{eqnarray*}
x_{y_{n+1}}
        & = & 2z_n+2^{n+1}-1\\
        & = & 2(n2^n)+2^{n+1}-1\\
        & = & n2^{n+1}+2^{n+1}-1\\
        & = & (n+1)2^{n+1}+1\\
        & = & z_{n+1}.
\end{eqnarray*}

Thus the induction hypothesis is true for $i$ and the lemma is proven.
\end{proof}

Now we use this information to estimate the growth rate of the
ball around 0.

\begin{lemma}
\label{right} For any $i>0$ and any $m$, $0 \leq m \leq 2^{i-1}$
we have that $d(0,m) \leq x_i$.
\end{lemma}

\begin{proof}
By induction we can see that the distance from any point to the
nearest vertex of level $E_k$ is at most $x_k$.  For any $m$ such
that $0 \leq m \leq 2^{i-1}$ the nearest vertex to $m$ of level
$k$ will lie in the interval $(0,2^{i-1})$. Thus
$$d(0,m) \leq min_k 2x_k + 2^{i-k-1}-1 =x_i.$$
\end{proof}
\begin{lemma}
\label{growth}
There is a function $G(j)$ defined below ($G(j)
\approx j^{.5\log j}$) such that
 $$G(j) \leq |B_{j}(0)| \leq 16 j^{2}G(j).$$
\end{lemma}

\begin{proof}
First for $i>2$ let
$$w_i=\sup_{k>0} 2^{k-1}+2^k((x_i-1)/2-x_k).$$
 (We want $i>2$ because all $x_i$ are odd except $x_2=2$.)
Thus $w_i$ is the largest integer such that there exists a path
from 0 to $w_i$ of length $(x_i-1)/2$.  This makes it is clear
that
$$B_{(x_i-1)/2}(0) \subset (-w_i,w_i).$$

Let $P$ be a path of length $(x_i-1)/2$ connecting 0 to $w_i$ and
$k$ be such that the longest step in $P$ is of size $2^k$. Suppose
that $w_i \geq 2^{i-2}$. The graph $G_k$ is symmetric about any
point of the form $w_i \pm l2^{k-1}$. Thus by combining $P$ and
the reflection of $P$ (about some suitably chosen point) we could
construct a path from 0 to $2^{i-1}$ of length at most $x_i-1$.
This is a contradiction. Thus
$$w_i<2^{i-2}$$
and
 \be B_{(x_i-1)/2}(0) \subset (-2^{i-2},2^{i-2}). \label{upper}
 \ee

On the other hand by Lemma \ref{right} gives us that
 \be
 [-2^{i-1},2^{i-1}] \subset B_{x_i}(0). \label{lower}
 \ee
Plugging $i=y_{n+1}$ into line (\ref{upper}) and $i=y_n$ into line
(\ref{lower}) gives
 \be
[-2^{y_n-1},2^{y_n-1}] \subset
 B_{x_{y_n}}(0) \subset B_{(x_{y_{n+1}}-1)/2}(0) \subset
(-2^{y_{n+1}-2},2^{y_{n+1}-2}).
 \ee
If $x_{y_{n}} \leq j < x_{y_{n+1}}$ then
$$B_{x_{y_n}}(0) \subset B_j(0) \subset B_{y_{n+1}}(0)$$ and
 \be \label{bounds}
[-2^{y_n-1},2^{y_n-1}] \subset
 B_{j} \subset
(-2^{y_{n+2}-2},2^{y_{n+2}-2}).
 \ee
We now rewrite this equation using the following definitions.  Let
$f(j)=\sup_n z_n \leq j$, $g(j)=2^{f(j)-1}$, $p(n)=(n^2+3n+2)/2$,
and $h(j)=2^{p(f(j)+2)-2}$.  Thus line (\ref{bounds}) becomes
 \be
[-g(j),g(j)] \subset B_{j} \subset (-h(j),h(j)).
 \ee

Then calculating
\begin{eqnarray*}
\frac{h(j)}{g(j)}&=&2^{p(f(j)+2)-2+p(f(j))-1}\\
    &=&2^{.5((f(j)^2+7f(j)+12)-(f(j)^2+3f(j)+2))-1}\\
    &=&2^{2f(j)+4}.
\end{eqnarray*}
By the definition of $z_n$ we get the bound
$$f(j)2^{2(j)}\leq j$$ and thus $f(j) < \log (j)$.  Putting these two
together we get that
$$\frac{h(j)}{g(j)}=2^{2f(j)+4}<2^{2\log(j)+4}=16j^2.$$
Thus
$$
[-g(j),g(j)] \subset
 B_{j} \subset
(-16j^2g(j)j^2,16j^2g(j)).
$$
Thus we can pick $G(j)=2g(j)+1.$  Finally we check that
$$G(j) \approx g(j) \approx 2^{.5f(j)^2} =j^{.5\log j}.$$

\end{proof}

\section{$\omega$-periodic graphs}
The graphs that we will consider are all obtain as the limit of
periodic graphs.
\begin{defn}
A graph $G$ with vertices labelled by $\Z$ is {\bf
$\omega$-periodic} if it is a union of periodic graphs over $\Z$.
\end{defn}

Our general result result about the growth of $\omega$-periodic
graphs is the following.
\begin{prop}\label{amenable}
Bounded degree $\omega$-periodic graphs are amenable.
\end{prop}
\begin{proof}
Let $G$ be an $\omega$-periodic graph with a uniformly bounded
degree. To show amenability it is enough to present a  growing
sets for which the ration between the size of the boundary of the
sets and the size of the sets is approaching $0$. $G$ is composed
of periodic layers ordered according to the density of the
vertices used. In particular the density of vertices that are
connected more than $k$ away, denoted by $k(s)$ exists and  is
going to $0$ with $k$. Hence if we consider a large interval of
size $n$, we get that it's boundary is smaller than $2k + k(s)n$.
Thus the ratio of the boundary to the interval can be made
arbitrarily small.
\end{proof}

\section{Polynomial Growth}

In this section we will use a subgraph $\tilde G$ of $G$ in
Section \ref{inter}. Again the vertices of $\tilde G$ are the
integers, $\Z$, and we define the sets of edges
$$E_0 =\{(i,i+1) : i \in \Z\}$$
and
$$
E_k=\left\{\left(2^k (n -1/2), 2^k(n+1/2)\right)\right\},
$$
for all $n \in \Z$ and $k >0$. We define the graph $\tilde G$ to
have edges
 $E= E_0 \bigcup \left( \cup_{k\geq 0} E_{2^k} \right).$

The proof that the volume of $\tilde B_m(0)$ grows polynomially in
$m$ is almost exactly like the proof of the volume growth of the
full graph in Section \ref{inter}.  First we calculate the
distance $\tilde x_{2^i}$ from $0$ to $2^{2^i-1}$.  Then we use
this information to bound the volume growth.  The difference is
that we get the formula
 $$\tilde x_{2^i} = \min_{0< k < i}2\tilde x_{2^k}+2^{2^i-2^k-1}-1.$$

We use the notation $\tilde B_j(0)$ to be the ball of radius $j$
in $\tilde G$ and
$$\tilde w_{2^i}
 =\max_k 2^{2^k-1} +
   2^{2^k}((\tilde x_{2^i}-1)/2- \tilde x_{2^k}).$$
 This gives us the following lemma.
\begin{lemma}\label{reader}
\begin{enumerate}
\item $\tilde x_{2^i}=2\tilde x_{2^{i-1}} + 2^{2^{i-1}-1} - 1$
\item $2^{2^{i-1}-1} \leq \tilde x_{2^i} \leq 2^{2^{i-1}}$ \label{r2}
\item $[-2^{2^i-1},2^{2^i-1}]\subset
    \tilde B_{\tilde x_{2^i}}(0),$ \label{r3}
\item $B_{2 \tilde x_{2^i}} \subset (-\tilde w_{2^i},\tilde
        w_{2^i})$ and \label{r4}
\item $\tilde w_{2^i} \leq 2^{2^i}$. \label{r5}
\end{enumerate}
\end{lemma}

\begin{proof}
The proof of these facts goes exactly as the proof of the
corresponding statements in Section \ref{inter}.
\end{proof}

\begin{lemma} If $j=\tilde x_{2^i}$ then
$$j^2 \leq |\tilde B_{j}(0)| \leq 8j^2.$$
\end{lemma}

\begin{proof}
The lower bound follows from condition \ref{r3} and the lower
bound in condition \ref{r2} in the previous lemma.  The upper
bound follows from conditions  \ref {r4}, \ref{r5} and the upper
bound in condition \ref{r2}.
\end{proof}

\section{Exponential Growth}

In this section we will construct an $\omega$-periodic graph that
contains a dyadic tree. Thus the graph has exponential volume
growth. Again we let
$$E_0 =\{(i,i+1) : i \in \Z\}$$
be the graph between adjacent integers.

Let $p_i$ be the $i$th prime,
  $$l_m= \prod_{i=1}^{i=2^m}(p_i)^{i}$$
  and
  $$t_{m,j}=(p_m)^j,\quad j=1\dots 2^{m-1}.$$
Notice that if $t_{m,j}=t_{m',j'}$ then $m=m'$ and $j=j'$.

 Define the $m$th level by
 $$E_m=\cup_{k\in\Z}\left(\cup_{j=1}^{2^{m-1}}
   \big((t_{m,j}+kl_{m+1},t_{m+1,2j-1}+kl_{m+1})
       \cup(t_{m,j}+kl_{m+1},t_{m+1,2j}+kl_{m+1})\big)\right).$$
Another way to describe $E_m$ is as follows.  Let
  $$V_m=\cup_{j=1}^{2^{m-1}} t_{m,j}.$$
Also let $L_m=V_m+\Z l_{m+1}$ and $R_m=V_{m+1}+\Z l_{m+1}$.
Then every edge
in $E_m$ has its leftmost endpoint in $L_m$ and its
rightmost endpoint in $R_m$.
Also every point in $L_m$ is the left hand end point of two edges in $E_m$.
Every point in $R_m$ is the right hand end point of one edge in $E_m$.

First we show that the graph contains a dyadic tree and then we
show that it has bounded degree.
\begin{lemma}
\label{tree}
There is a dyadic tree rooted at $t_{1,1}=2.$
\end{lemma}

\proof
Note that $V_{m+1} \subset L_{m+1} \cap R_m$.  Then the $2^m$
vertices at distance $m$ from the root are $V_{m+1}$.
 \qed

\begin{lemma}
\label{disjoint}
 $L_m \cap \left(\cup_{j=1}^{m-1} L_j\right)=\emptyset$ and
 $R_m \cap \left(\cup_{j=1}^{m-1} R_j\right)=\emptyset$.
\end{lemma}

\proof
 Fix an $m$.  Every element of  $L_m \mod l_{m+1}$ is only
divisible by powers of $p_m$.  Every element of
 $\cup_{j=1}^{m-1} L_{j} \mod l_{m+1}$
is divisible by at least one prime less than or equal to $p_{m-1}$. Thus
the two sets are disjoint.

Fix an $m$.  Every element of  $R_m \mod l_{m+1}$ is only
divisible by powers of $p_{m+1}$.  Every element of
 $\cup_{j=1}^{m-1} R_{j} \mod l_{m+1}$
is divisible by at least one prime less than or equal to $p_{m}$. Thus the
two sets are disjoint.\qed

\begin{lemma}
\label{five}
The degree of any vertex in $E$ is at most five.
\end{lemma}

\proof
For any $z \in \Z$ the degree of $z$ is 2 plus twice the number of
$m$ such that $z \in L_m$ plus the number of $m$ such that $z \in
R_m$. Thus by Lemma \ref{disjoint} the degree of a vertex is at
most five. \qed

\end{document}